\documentclass[reqno, 11pt]{amsart}
\usepackage[utf8]{inputenc}
\usepackage[T1]{fontenc}
\usepackage{lmodern}
\usepackage{epsfig,amsmath, amsthm ,amsfonts,latexsym,stmaryrd, amssymb}
\usepackage[abs]{overpic}
\usepackage[dvipsnames]{xcolor}
\usepackage[linktocpage=true]{hyperref}
\usepackage{caption}
\usepackage{subcaption}
\usepackage{dsfont}
\usepackage{mathtools}
\usepackage{tikz-cd}
\usetikzlibrary{shapes,arrows}
\usepackage{comment}
\usepackage{enumitem}
\usepackage{marginnote}
\usepackage[all,cmtip]{xy}
\usepackage[backend=biber, maxnames=4, style=alphabetic]{biblatex}
\hypersetup{
  colorlinks,
  citecolor=teal,
  linkcolor=purple,
  urlcolor=teal}
\theoremstyle{plain}

\newtheorem{dummy}{anything}[section]
\newtheorem{thm}{Theorem}

\theoremstyle{definition}
\newtheorem{definition}[dummy]{Definition}

\newtheorem{remark}[dummy]{Remark}

\theoremstyle{remark}
\newcommand{\R}{\mathbb{R}}

\def\R{\mathbb{R}}

\DeclareFontFamily{U}{mathx}{}
\DeclareFontShape{U}{mathx}{m}{n}{<-> mathx10}{}
\DeclareSymbolFont{mathx}{U}{mathx}{m}{n}
\DeclareMathAccent{\widecheck}{0}{mathx}{"71}

\begin{document}

\author{Agustin Moreno}

\address[A.\ Moreno]{Universit\"at Heidelberg \\ Mathematisches Institut \\ Heidelberg \\ Germany}

\email{\href{mailto:agustin.moreno2191@gmail.com}{agustin.moreno2191@gmail.com}}

\date{\today}



\title{Correspondence between orbit equivalences and exact symplectomorphisms for Anosov flows}

\begin{abstract}
    We show that for a closed hyperbolic $3$-manifold, two Anosov flows are orbital equivalent if and only if the corresponding Anosov Liouville structures are exact symplectomorphic.
\end{abstract}

\maketitle

\tableofcontents

\section{Introduction}

In this note, we establish a two-way relationship between Anosov flows on a closed hyperbolic $3$-manifold $M$, and that of certain symplectic structures on $\mathbb R \times M$. We first need the following notion.

\begin{definition} \label{def:alstructure}
An \emph{Anosov Liouville structure} on $\R \times M$ is a Liouville structure induced by a smooth Liouville form $\lambda$ of the form
\begin{align} \label{eq:exponentialanosov}
\lambda = e^{-s} \alpha_- + e^s \alpha_+
\end{align}
where $(\alpha_-, \alpha_+)$ is a pair of contact forms on $M$ satisfying the following two conditions.
\begin{enumerate}
    \item The contact plane fields $\xi_\pm = \ker \alpha_\pm$ are transverse everywhere.
    \item The line distribution $\xi_- \cap \xi_+$ is spanned by an Anosov vector field.
\end{enumerate}
\end{definition}

The pair $(\xi_-,\xi_+)$ constitutes a \emph{bi-contact structure} on $M$, i.e., a
transverse pair of contact plane fields with opposite orientations, which is said to \emph{support} the corresponding Anosov flow. Every Anosov vector field can be presented in this way (see Mitsumatsu \cite[Theorem 3]{M95}, Eliashberg--Thurston~\cite{ET98}, Hozoori~\cite{H22a} and Massoni ~\cite{Mas22}).

The natural notion of equivalence between two Anosov flows is that of \emph{orbital equivalence}, i.e.\ the existence of a \emph{homeomorphism} taking orbits to orbits. On the symplectic side, the notion of equivalence between two Anosov Liouville structures $\lambda,\lambda'$ is that of \emph{exact symplectomorphism}, i.e.\ the existence of a \emph{diffeomorphism} $f$ satisfying $f^*\lambda = \lambda' +dh$, where $h$ is a compactly supported smooth function. 

\section{Main theorem}

In what follows, we will implicitly assume that all Anosov flows have orientable weak-stable and weak-unstable foliations. This can be achieved by passing to a double cover. The goal of this note is to prove the following.

\begin{thm}
    Let $M$ be a closed hyperbolic $3$-manifold. Then there is a $1:1$ correspondence
    $$
    \frac{\{\mbox{Anosov flows}\}}{\mbox{orbital equivalence}}\longleftrightarrow \frac{\{\mbox{Anosov Liouville structures}\}}{\mbox{exact symplectomorphism}}
    $$
\end{thm}

In other words, counting orbital equivalence classes of Anosov flows is the same as counting exact sympleectomorphism classes of Anosov Liouville structures (which in the hyperbolic case are finite by the main result of Bowden--Colin \cite{BC26} and independently Gabai--Li \cite{GL26}). The direct arrow was obtained in recent work of Bowden--Massoni \cite{BM26}, so it suffices to prove the other arrow, i.e.\ that exact symplectomorphic Anosov Liouville structures induce orbital equivalent Anosov flows. This answers Question 5.4 of Massoni's survey \cite{Mas25}, in the hyperbolic case.

\begin{proof}
Let $f:\mathbb{R} \times M' \rightarrow \mathbb{R} \times M$ be an exact symplectomorphism, i.e.\ satisfying $f^*\lambda=\lambda' + dh$ for some compactly supported smooth function $h:M'\rightarrow \mathbb{R}$. Let $\xi_\pm,\xi_\pm'$ be the corresponding contact structures at the ends.

Since $h$ vanishes outside of a compact set, $f$ induces contactomorphisms
$$
g_+ : (M',\xi_+')\rightarrow (M,\xi_+),\; g_-: (M',\xi_-')\rightarrow (M,\xi_-)
$$
Moreover, $g_\pm$ arise by restricting the same diffeomorphism $f$ to slices, so they are homotopic. By Gabai–Meyerhoff–Thurston \cite{GMT}, homotopic diffeomorphisms on a hyperbolic $3$-manifold are isotopic, so if $g:=g_+$, then $g^*\xi_\pm$ is isotopic to $\xi_\pm'$.

By a recent result of Bowden--Colin (Theorem 1.2 in \cite{BC26}), the corresponding flows are (isotopy) orbit equivalent. This finishes the proof.
\end{proof}

\begin{remark}
   By work of Barthelm\'e--Frankel--Mann \cite{BFM}, the orbital equivalence class of an Anosov flow $\psi$ on a hyperbolic $3$-manifold $M$ is completely determined by the set $\mathcal{P}(\psi)$ of unoriented free homotopy classes of loops in $M$ represented by periodic orbits (the \emph{free homotopy data}). Our main result then implies that, for hyperbolic $M$, two Anosov flows $\psi_1,\psi_2$ admit an automorphism $\Phi$ of $\pi_1(M)$ such that $\Phi_*\mathcal{P}(\psi_1)=\mathcal{P}(\psi_2)$ if and only if the corresponding Anosov Liouville structures are exact symplectomorphic. In other words, for a closed hyperbolic $3$-manifold, there is a $1:1$ correspondence
   $$
    \frac{\{\mbox{free homotopy data}\}}{\mbox{bijection}}\longleftrightarrow \frac{\{\mbox{Anosov flows}\}}{\mbox{orbital equivalence}}\longleftrightarrow \frac{\{\mbox{Anosov Liouville structures}\}}{\mbox{exact symplectomorphism}}.
    $$
\end{remark}

\textbf{Acknowledgments.} The author has received funding of the Air Force Office of Scientific Research under award number FA8655-24-1-7012.

\printbibliography

\end{document}